\documentstyle[11pt]{article}
\def\Empty{}
\catcode`\@=11

\def\section{\@startsection {section}{1}{\z@}{-3.5ex plus -1ex minus
-.2ex}{2.3ex plus .2ex}{\large\bf}}

\def\eqalign#1{\,\vcenter{\openup\jot\m@th
  \ialign{\strut\hfil$\displaystyle{##}$&$\displaystyle{{}##}$\hfil
        \crcr#1\crcr}}\,}

\def\mydesc{\list{}{\labelwidth\z@ \itemindent-\leftmargin
\listparindent 1.5em
\let\makelabel\descriptionlabel}}

\def\fnum@figure{{\small Figure \thefigure}}
\def\fakefigure{\def\@captype{figure}}
\long\def\@makecaption#1#2{
    \vskip 10pt
    \def\FCap{#2} \def\NoCap{\ignorespaces}
    \ifx \FCap\NoCap
       \setbox\@tempboxa\hbox{#1}  
      \else
       \setbox\@tempboxa\hbox{#1: \small \it #2}
    \fi
    \ifdim \wd\@tempboxa >\hsize   
        \unhbox\@tempboxa\par      
      \else                        
        \hbox to\hsize{\hfil\box\@tempboxa\hfil}
    \fi}

\def\@oddhead{\hbox{}\rightmark \hfil \rm\thepage}
\def\sectionmark#1{\markright {\sc{\ifnum \c@secnumdepth >\z@
      \S\thesection.\hskip 1em\relax \fi #1}}}

\catcode`\@=12

\def\oplabel#1{
  \def\OpArg{#1} \ifx \OpArg\Empty {} \else
        \label{#1}
  \fi}
\def\MakeStEnv#1{
  \newenvironment{#1}[2]{
  \begin{#1St} \oplabel{##1}%
  \global\def\CrntSt{\thetheoremSt}%
  {\rm ##2}%
}{
  \end{#1St} }
}
\MakeStEnv{theorem}
\MakeStEnv{corollary}
\MakeStEnv{proposition}
\MakeStEnv{lemma}
\MakeStEnv{fact}
\MakeStEnv{conjecture}
\MakeStEnv{comment}
\MakeStEnv{remark}
\newenvironment{proof}[1]{
  \def\PfArg{#1}
  \ifx\PfArg\Empty
        \edef\PfArg{\CrntSt}  \fi
 \startproof{\PfArg}%
}{
  \finishproof{\PfArg}
}
\newcommand{\startproof}[1]{
  \medbreak\mbox{}
  {\it Proof of #1:}%
}
\newcommand{\finishproof}[1]{
  \def\FPArg{#1}
  \ifx\FPArg\Empty
        \def\FPArg{\CrntSt}  \fi
  \smallbreak\noindent\makebox[\textwidth]{\hfill\fbox{\FPArg}}
  \medbreak\noindent
}
\def\hyperbolic{{\bf H}}

\def\H{\hyperbolic}

\def\PSL#1{{\rm PSL}_{#1}}

\def\til{\widetilde}
\def\hat{\widehat}

\def\st{{\rm st}}

\def\fix{{\rm fix}}
\def\inj{{\rm inj}}

\def\PSL{{\rm PSL}}

\def\rmd{{\rm d}}

\def\bfC{{\bf C}}

\def\bfH{{\bf H}}

\def\bfZ{{\bf Z}}

\def\calD{{\cal D}}
\def\calH{{\cal H}}

\def\rs{\overline{\bfC}}

\title{Cores of hyperbolic $3$-manifolds and limits of Kleinian groups
II}
\author{James W. Anderson and Richard D. Canary\thanks{Research
supported in part by the National Science Foundation and a fellowship
from the Sloan Foundation}}
\date{\small{Faculty of Mathematical Studies, University of
Southampton\\ Highfield, Southampton, SO17 1BJ, England}\\
\small{Department of Mathematics, University of Michigan, Ann Arbor,
MI 48109}}
\begin{document}

\maketitle

\section{Introduction}
\label{introduction}

Troels J\" orgensen conjectured that the algebraic and geometric
limits of an algebraically convergent sequence of
isomorphic Kleinian groups agree if there are {\em no new parabolics}
in the algebraic limit.  We prove that this conjecture holds in 
``most'' cases.  In particular,  we show that it holds when the domain of
discontinuity of the algebraic limit of such a sequence is non-empty
(see Theorem \ref{strong-limit}).
We further show, with the same assumptions,
that the limit sets of the groups in the sequence
converge to the limit 
set of the algebraic limit.  As a corollary, we verify the conjecture
for finitely generated Kleinian groups which are not (non-trivial)
free products of surface groups and infinite cyclic groups
(see Corollary \ref{group-condition}).  These
results are extensions of similar results for purely loxodromic groups
which can be found in \cite{anderson-canary}.  Thurston
\cite{thurston-notes} previously established these results in the case
that the Kleinian groups are freely indecomposable (see also Ohshika
\cite{ohshika-caratheodory}, \cite{ohshika-strong},
\cite{ohshika-divergent}).  Using different techniques than ours,
Ohshika \cite{ohshika-limits} has proven versions of these results for purely
loxodromic function groups. 

\medskip

Both authors would like to thank the Institut Henri Poincar\'e for its
hospitality during the writing of this paper, as well as the referee
for useful comments.

\section{Preliminaries}

The purpose of this section is to describe the background material
used in this paper.

\subsection{Convergence of Kleinian groups}

A {\em Kleinian group} is a discrete subgroup of $\PSL_2(\bfC)$, which
we view as acting either on hyperbolic $3$-space $\bfH^3$ via
isometries or on the Riemann sphere $\rs$ via M\"obius transformations.
The action of $\Gamma$ partitions $\rs$ into the {\em domain of
discontinuity} $\Omega(\Gamma)$, which is the largest open subset of
$\rs$ on which $\Gamma$ acts properly discontinuously, and the {\em
limit set} $\Lambda(\Gamma)$.  A Kleinian group is {\em
non-elementary} if its limit set contains at least three points, and
is {\em elementary} otherwise.  A Kleinian group is elementary if and
only if it is virtually abelian; recall that a group is {\em virtually
abelian} if it contains a finite index abelian subgroup.  We refer
the reader to Maskit \cite{maskit-book} for a more detailed discussion
of the theory of Kleinian groups.

It is often convenient to view a Kleinian group as the image of a
discrete, faithful representation of a group into
$\PSL_2(\bfC)$. Given a finitely generated group $G$, let $\calD(G)$
denote the space of all discrete, faithful representations of $G$ into
$\PSL_2(\bfC)$.  A sequence $\{\rho_j\}$ in $\calD(G)$ converges {\em
algebraically} to $\rho$ if $\{\rho_j(g)\}$ converges to $\rho(g)$ for
each $g\in G$.  It is a fundamental result of J\o rgensen
\cite{jorgensen} that $\calD(G)$ is a closed subset
of $Hom(G,\PSL_2 (\bfC))$ when $G$ is finitely generated and not
virtually abelian. 

An isomorphism $\alpha: \Phi\rightarrow\Theta$ between Kleinian
groups $\Phi$ and $\Theta$ is {\em type-preserving} if $\varphi$ is
parabolic if and only if $\alpha(\varphi)$ is parabolic for all
$\varphi\in\Phi$.  More generally, an algebraically convergent
sequence $\{\rho_j\}\subset\calD(G)$ with limit $\rho\in\calD(G)$ is
{\em type-preserving} if $\rho_j\circ\rho^{-1}:
\rho(G)\rightarrow\rho_j(G)$ is type-preserving for all $j$.

A sequence $\{\Gamma_j\}$ of Kleinian groups converges {\em geometrically}
to a Kleinian group $\hat\Gamma$ if every element of
$\hat\Gamma$ is the limit of a sequence $\{\gamma_j\in\Gamma_j\}$ and
if every accumulation point of every sequence
$\{\gamma_j\in\Gamma_j\}$ lies in $\hat\Gamma$.  We make use of the
following proposition, which assures that algebraically convergent
sequences in ${\cal D}(G)$ have geometrically convergent subsequences.

\begin{proposition}{alg-geom}{(Proposition 3.8 in J\o rgensen and Marden
\cite{jorgensen-marden})}{}
If $G$ is not virtually abelian and
$\{\rho_j\}\subset \calD(G)$ is an algebraically convergent sequence
with limit $\rho$, then there exists a subsequence $\{\rho_{j_k}\}$ of
$\{\rho_j\}$ so that $\{\rho_{j_k}(G)\}$ converges geometrically to a
Kleinian group $\hat\Gamma$ with $\rho(G)\subset\hat\Gamma$.
\end{proposition}

A sequence of closed sets $\{ X_j\}$ in $\rs$ {\em converges in the
Hausdorff topology} to a closed set $X$ in $\rs$ if every point of $X$
is the limit of a sequence of points $\{x_j\in X_j\}$ and if every
accumulation point of a sequence $\{x_j\in X_j\}$ is contained in
$X$.  With this topology, the collection of non-empty closed subsets
of $\rs$ is compact.  Convergence of closed sets is always assumed to
be in the Hausdorff topology.

If $\{\rho_j\}\subset\calD$ converges algebraically to $\rho$ and
$\{\rho_j(G)\}$ converges geometrically to $\rho(G)$, we say that
$\{\rho_j\}$ converges {\em strongly} to $\rho$.
It is conjectured that if $G$ is not virtually abelian and if
$\{\rho_j\}\subset\calD(G)$ converges algebraically to $\rho$, then
$\{\Lambda(\rho_j(G))\}$ converges to $\Lambda(\rho(G))$ if and only
if $\{\rho_j\}$ converges strongly to $\rho$.  We make use of the
following partial result in the direction of this conjecture.

\begin{proposition}{limcon}{(Proposition 4.2 in J\o rgensen and Marden
\cite{jorgensen-marden})}{} Let $G$ be a finitely generated group which
is not virtually abelian.
If $\{\rho_j\}\subset {\cal D}(G)$ converges algebraically to
$\rho\in{\cal D}(G)$, $\Omega(\rho(G))$ is non-empty, and
$\{\Lambda(\rho_j(G))\}$ converges to $\Lambda(\rho(G))$, then
$\{\rho_j\}$ converges strongly to $\rho$. 
\end{proposition}

The following lemma indicates the geometric significance, on the level
of the quotient manifolds, of the geometric  convergence of a sequence
of Kleinian groups.  (For a proof, see Theorem 3.2.9 of
Canary, Epstein, and Green \cite{canary-epstein-green}, and Theorem
E.1.13 and Remark E.1.19 of Benedetti and Petronio
\cite{benedetti-petronio}.)  Let $0$ denote a fixed choice of
basepoint for $\H^3$, and let $B_{R}(0)$ denote the ball of radius $R$
centered at $0$.


\begin{lemma}{fi convergence}{}{} A sequence of torsion-free Kleinian
groups $\{\Gamma_j\}$ converges geometrically to a torsion-free
Kleinian group $\hat\Gamma$ if and only if there exists a sequence
$\{(R_j,K_j)\}$ and a sequence of maps $\til f_j:B_{R_j}(0) \to \H^3$
such that 
\begin{enumerate}
\item $R_j\to\infty$ and $K_j\to 1$ as $j\rightarrow\infty$;
\item the map $\til f_j$ is a $K_j$-bilipschitz diffeomorphism onto
its image, $\til f_j(0)=0$, and for any compact subset $A$ of $\H^3$,
$\{\til f_j|_A\}$ converges to the identity; 
\item if $V_j =B_{R_j}(0)/\Gamma_j$, then $V_j$ is a submanifold of
$N_j =\bfH^3/\Gamma_j$ and $\til f_j$ descends to a map $f_j:V_j \to
\hat N$, where $\hat N = \H^3/\hat\Gamma$; moreover, $f_j$ is also a
$K_j$-bilipschitz diffeomorphism onto its image.
\end{enumerate}
\end{lemma}

\subsection{Types of Kleinian groups}

Given a set $X\subset\bfH^3\cup\rs$ and a Kleinian group $\Gamma$,
define the {\em stabilizer} of $X$ in $\Gamma$ to be
\[ \st_{\Gamma}(X) =\{ \gamma\in\Gamma\: :\: \gamma(X) =X\}. \]
A {\em component subgroup} of $\Gamma$ is the stabilizer in $\Gamma$
of a component $\Delta$ of $\Omega(\Gamma)$.  A set
$X\subset\bfH^3\cup\rs$ is {\em precisely invariant} under a subgroup
$\Phi$ of $\Gamma$ if $\st_{\Gamma}(X) =\Phi$ and if $X\cap\gamma(X)$
is empty for all $\gamma\in\Gamma -\Phi$.

There are several classes of Kleinian groups of particular interest in
our paper.  A {\em quasifuchsian} group is a finitely generated
Kleinian group whose limit set is a Jordan curve and which contains
no element interchanging the components of its domain of
discontinuity, while an {\em extended quasifuchsian} group is a
finitely generated Kleinian group whose limit set is a Jordan curve
and which does contain an element interchanging the components of its
domain of discontinuity.  Note that an extended quasifuchsian group
contains a canonical quasifuchsian subgroup of index two, namely the
component subgroup associated to either of the components of its
domain of discontinuity.

A {\em degenerate} group is a finitely generated Kleinian group whose
domain of discontinuity and limit set are both non-empty and
connected.  A {\em web} group is a finitely generated Kleinian group
whose domain of discontinuity contains infinitely many components,
and each component subgroup is quasifuchsian.  A finitely
generated Kleinian group is a {\em generalized web} group if it is
either quasifuchsian, extended quasifuchsian or a web group.


The convex core $C(N)$ of a hyperbolic $3$-manifold $N=\H^3/\Gamma$ is
the quotient of the convex hull $CH(\Lambda(\Gamma))$ of the limit set
of $\Gamma$ by $\Gamma$. (The convex core can also be defined to be
the smallest convex subset of $N$ whose inclusion is a homotopy
equivalence.)  A finitely generated  Kleinian group is {\em
geometrically finite} if its convex core has finite volume.  A
torsion-free Kleinian group is {\em topologically tame} if its
quotient $3$-manifold $N =\bfH^3/\Gamma$ is homeomorphic to the
interior of a compact $3$-manifold.  We note that geometrically finite
Kleinian groups are topologically tame (see Marden \cite{marden}).  A
theorem of Scott \cite{scott} guarantees that any hyperbolic
3-manifold with finitely generated fundamental group contains a
compact submanifold, called a {\em compact core}, whose inclusion is a
homotopy equivalence. It is conjectured (see Marden \cite{marden})
that every hyperbolic 3-manifold with finitely generated fundamental
group is homeomorphic to the interior of its compact core, and hence
topologically tame.

\subsection{Decompositions of Kleinian groups}
\label{decomp}

In Section \ref{strong limits}, we make use of two related
decompositions, due to Abikoff and Maskit \cite{abikoff-maskit}, of a
finitely generated, torsion-free Kleinian group.
We first discuss how 
a non-elementary, finitely generated, torsion-free Kleinian  group with
connected limit set and non-empty domain of discontinuity can be built from
generalized web groups and degenerate groups without accidental parabolic
elements.  We then discuss the
decomposition of a finitely generated Kleinian group into groups with
connected limit sets and elementary groups.  

We begin with a few definitions.  Let $\Gamma$ be a torsion-free
Kleinian group, and let $\gamma\in\Gamma$ be a parabolic element.  A
{\em cusp region} for $\gamma$ is a closed Jordan domain $D$ which is
precisely invariant under $\langle\gamma\rangle$ in $\Gamma$ and
intersects $\Lambda(\Gamma)$ only at the fixed point of $\gamma$.
Note that this immediately implies that $\gamma$ is primitive in
$\Gamma$ and that $\gamma$ cannot lie in a rank two abelian subgroup
of $\Gamma$.  The interior of a cusp region descends to a punctured
disc neighborhood of a cusp on the Riemann surface
$\Omega(\Gamma)/\Gamma$. 

As the choice of a cusp region is by no means canonical, we introduce
a notion of equivalence.  Say that cusp regions $D$ and $D'$ for
parabolic elements $\gamma$ and $\gamma'$ of $\Gamma$ are {\em
equivalent} if their images on $\Omega(\Gamma)/\Gamma$ are
neighborhoods of the same cusp on $\Omega(\Gamma)/\Gamma$.  This
implies in particular that $\gamma$ and $\gamma'$ are conjugate in
$\Gamma$; we note that the converse need not hold, as a single
primitive parabolic element may stabilize two inequivalent cusp
regions.

Let $\Gamma$ be a non-elementary, finitely generated, torsion-free
Kleinian group with connected limit set; in particular, every
component of $\Omega(\Gamma)$ is simply connected.  Let $\Delta$ be a
component of $\Omega(\Gamma)$, and let $\Phi =\st_{\Gamma}(\Delta)$ be
its component subgroup.  By the classical Uniformization Theorem for
Riemann surfaces, there exists a conformal homeomorphism $f:
\Delta\rightarrow\bfH^2$.  The elements of $f \Phi f^{-1}$ are 
conformal homeomorphisms of $\bfH^2$, and so are necessarily M\"obius
transformations.  An {\em accidental parabolic element} $\theta$ of
$\Phi$ is a parabolic element whose conjugate $f\theta f^{-1}$ by $f$
is a primitive hyperbolic element of $f\Phi f^{-1}$. Denoting by
$l$ the line in $\bfH^2$ joining the fixed points of $f\theta
f^{-1}$, we define the {\em axis} of $\theta$ to be $c_{\theta}
=f^{-1}(l)$ and the {\em full axis} of $\theta$ to be the Jordan curve
$C_{\theta} =c_{\theta}\cup\fix(\theta)$; for a detailed discussion of
accidental parabolics, we refer the reader to Maskit
\cite{maskit-book}, particularly Chapter IX.D.10. By construction,
$c_{\theta}$ is precisely invariant under $\Theta
=\langle\theta\rangle$ in $\Gamma$, and $C_{\theta}$ separates
$\Lambda(\Phi)$, and hence separates $\Lambda(\Gamma)$.

Before describing the decomposition of a Kleinian group along an
accidental parabolic element, we make a couple of remarks.  First, a
given parabolic element of $\Gamma$ may lie in several component
subgroups, being accidental in some and not accidental in others; 
hence, we make the convention that whenever we refer to an accidental
parabolic element of $\Gamma$, we actually refer to the parabolic
element along with the implicit choice for the component of
$\Omega(\Gamma)$ containing its axis.  Second, a primitive parabolic
element of a Kleinian group which keeps invariant a component of the
domain of discontinuity must either be accidental in that component or
must have a cusp region in that component.  Third, a finitely
generated, non-elementary Kleinian group with connected limit set and
non-empty domain of discontinuity contains no accidental parabolic
elements if and only if it is either a degenerate group without accidental
parabolics or a generalized web group (this follows immediately from
Theorem 4 of Maskit \cite{maskit-boundaries}).

\medskip

Given an accidental parabolic element $\theta$ of $\Gamma$ with full
axis $C_{\theta}$, let $P_1$ and $P_2$ be the two components of $\rs
-\Gamma(C_{\theta})$ which contain $C_{\theta}$ in their boundaries,
and let $\Gamma_m$ be the stabilizer of $P_m$ in $\Gamma$.  Let $E_1$
and $E_2$ be the closed discs in $\rs$ determined by $C_{\theta}$,
labeled so that $C_{\theta}$ separates the interior of $E_m$ from
$P_m$.  

Suppose that $P_1$ and $P_2$ are not equivalent under $\Gamma$.
In this case, we say that $C_{\theta}$ is a {\em separating} full axis
for $\theta$, and we refer to $\Gamma_1$ and $\Gamma_2$ as the {\em
factor subgroups} of this decomposition.  It is easy to show that
$\Gamma$ is generated by its factor subgroups $\Gamma_1$ and
$\Gamma_2$, that both factor subgroups are non-elementary and finitely
generated, that $\Gamma_1\cap\Gamma_2 =\langle\theta\rangle$, and that
$C_{\theta}$ is precisely invariant under $\langle\theta\rangle$ in
$\Gamma_m$ for both $m$.  The statement of the Klein-Maskit
combination theorem of type I given below, adapted from the statement
in Maskit \cite{maskit-book}, essentially states that this operation
of decomposition can be reversed.

Recall that a {\em fundamental domain} for a finitely generated
Kleinian group $\Gamma$ is an open subset $D$ of $\Omega(\Gamma)$ so
that $D$ is precisely invariant under the identity, every point of
$\Omega(\Gamma)$ is equivalent under the action of $\Gamma$ to a point
of $\overline{D}$ and the boundary of $D$ is a finite collection of
analytic arcs.

\begin{theorem}{kmcombI}{}{(Klein-Maskit combination I)}  Let
$\Gamma_1$ and $\Gamma_2$ be non-elementary, finitely generated
Kleinian groups whose intersection $\Gamma_1\cap\Gamma_2$ is the
parabolic cyclic group $\langle\theta\rangle$.  Let $C$ be a Jordan
curve in $\rs$ determining closed discs $E_1$ and $E_2$, so that $E_m$
is a cusp region for $\theta$ in $\Gamma_m$.

Then, the following hold.  First, the group
$\Gamma =\langle\Gamma_1,\Gamma_2\rangle$ is a Kleinian group
isomorphic to the amalgamated free product of $\Gamma_1$ and
$\Gamma_2$ along $\langle\theta\rangle$.  Second, if $D_m$ is a
fundamental domain for $\Gamma_m$ so that $D_m\cap E_m$ is a
fundamental domain for the action of $\langle\theta\rangle$ on $E_m$,
$D_m\cap E_{3-m}$ has non-empty interior, and $D_1\cap C_{\theta} =
D_2\cap C_{\theta}$, then $D = (D_1\cap E_2)\cup (D_2\cap E_1)$ is a
fundamental domain for the action of $\Gamma$.  Third, every cusp
region for $\Gamma_m$ which does not intersect any
$\Gamma_m$-translate of $E_m$ is a cusp region for $\Gamma$, and every
cusp region for $\Gamma$ is a cusp region for either
$\Gamma_1$ or $\Gamma_2$.
\end{theorem}

Suppose, on the other hand, that $P_1$ and $P_2$ are equivalent by
$\xi\in\Gamma$, so that $\xi(P_1) =P_2$.  In this case we say that
$C_{\theta}$ is a {\em non-separating} full axis for $\theta$, and we
refer to $\Gamma_2$ and $\langle\xi\rangle$ as the {\em factor
subgroups} of this decomposition.  It is easy to show that $\Gamma$ is
generated by its factor subgroups $\Gamma_2$ and $\langle\xi\rangle$,
that $\Gamma_2$ is finitely generated and non-elementary, that
$\Gamma_2\cap\xi^{-1}\Gamma_2\xi =\langle\theta\rangle$, and that
$C_{\theta}$ is precisely invariant under $\langle\theta\rangle$ in
$\Gamma_2$ (even though we may have that $\xi(C_{\theta})\cap
C_{\theta}$ is non-empty if $\xi$ and $\theta$ commute).  The
statement of the Klein-Maskit combination theorem of type II, adapted
from the statement in Maskit\cite{maskit-book}, essentially states
that this decomposition can be reversed.

\begin{theorem}{kmcombII}{}{(Klein-Maskit combination II)}  Let
$\Gamma^0$ be a non-elementary, finitely generated Kleinian group, and
let $\langle\theta_1\rangle$ and $\langle\theta_2\rangle$ be
parabolic cyclic subgroups of $\Gamma^0$.  Let $C_1$ and $C_2$ be
Jordan curves in $\rs$ so that, for both $m$, $C_m$ determines a
closed disc $E_m$ which is a cusp
region for $\langle\theta_m\rangle$ in $\Gamma^0$.  We also require
that $E_1$ and $E_2$ are inequivalent cusp regions for $\Gamma^0$.  Let
$\xi$ be a M\"obius transformation so that $\xi\theta_1\xi^{-1}
=\theta_2$, $\xi(C_1) =C_2$, and the image of the interior of $E_1$
under $\xi$ is disjoint from the interior of $E_2$.

Then, the following hold.  First, the group
$\langle\Gamma^0,\xi\rangle$ is a Kleinian group isomorphic to the HNN
extension of $\Gamma^0$ by $\xi$.  Second, if $D^0$ is a fundamental 
domain for $\Gamma^0$ so that $D^0\cap E_m$ is a fundamental domain
for the action of $\langle\theta_m\rangle$ on $E_m$ and $\xi(D_1\cap
C_1) =D_2\cap C_2$, then $D =D^0\cap A$ is a fundamental domain for
$\Gamma$, where $A =\rs -(E_1\cup E_2)$.  Third, every cusp region for
$\Gamma^0$ which does not intersect any $\Gamma^0$-translate of $E_1$
or $E_2$ is a cusp region for $\Gamma$, and every cusp region for
$\Gamma$ is a cusp region for $\Gamma^0$.
\end{theorem}

Abikoff and Maskit \cite{abikoff-maskit} showed that given a finitely
generated, non-elementary Kleinian group $\Gamma$ with connected limit
set and non-empty domain of discontinuity, one can produce a finite
collection $\{\Gamma_1,\ldots,\Gamma_k,\Gamma_{k+1}, \ldots,
\Gamma_l\}$ of subgroups of $\Gamma$, where $\Gamma_i$
is either a generalized web group or a degenerate group without
accidental parabolic elements for $i\le k$, and $\Gamma_i$ is a cyclic
group for $i>k$, so that $\Gamma$ can be built from
$\{\Gamma_1,\ldots,\Gamma_l\}$ by repeatedly performing 
Klein-Maskit combinations of types I and II.  Specifically, set
$\Gamma^1=\Gamma_1$; for $j\le k$, form $\Gamma^j$ from $\Gamma^{j-1}$
and $\Gamma_j$ by using a Klein-Maskit combination of type I along 
a common parabolic cyclic subgroup of $\Gamma^{j-1}$ and $\Gamma_j$,
and for $j>k$, form $\Gamma^j$ from $\Gamma^{j-1}$ and $\Gamma_j$ by a
Klein-Maskit combination of type II, where the generator of $\Gamma_j$
pairs inequivalent cusp regions of $\Gamma^{j-1}$.  The final result
of this process $\Gamma^l$ is the original group $\Gamma$.

\medskip

Abikoff and Maskit \cite{abikoff-maskit} also showed that any
torsion-free, finitely generated Kleinian group with non-empty domain
of discontinuity can be constructed from a finite collection of
elementary groups and groups with connected limit set by a finite
number of applications of the Klein combination theorem.  We recall
the statement of the Klein combination theorem below for reference,
and refer the reader to Maskit \cite{maskit-book} for a complete
discussion.

\begin{theorem}{theo-klein-combination}{(Klein combination theorem)}{}
Let $\Gamma_1$ and $\Gamma_2$ be Kleinian groups with non-empty
domains of discontinuity, and suppose there exist fundamental domains
$D_1$ and $D_2$ for $\Gamma_1$ and $\Gamma_2$ so that each contains
the exterior of the other.  Then, $\Gamma
=\langle\Gamma_1,\Gamma_2\rangle$ is a Kleinian group isomorphic to
the free product $\Gamma_1\ast\Gamma_2$, and $D =D_1\cap D_2$ is a
fundamental domain for $\Gamma$.
\end{theorem}

\subsection{Relative compact cores and ends of Kleinian groups}

A {\em horoball} in $\bfH^3$ is an open Euclidean ball (or half-space)
$B$ in $\bfH^3$ whose (Euclidean) 
closure in $\bfH^3\cup\rs$ intersects $\rs$ in a single
point, which is the {\em center} of the horoball.  A {\em precisely
invariant system of horoballs} ${\cal H}$ for a Kleinian group
$\Gamma$ is a $\Gamma$-invariant collection of disjoint horoballs
centered at parabolic fixed points of $\Gamma$, such that there is one
based at every parabolic fixed point. It is a consequence of the
Margulis Lemma (see Benedetti and Petronio \cite{benedetti-petronio}
or Maskit \cite{maskit-book}) that every Kleinian group has a
precisely invariant system of horoballs.  Given a precisely invariant
system of horoballs $\calH$ for $\Gamma$, set $N_{\cal H} =(\H^3-{\cal
H})/\Gamma$.

A {\em relative compact core} of $N_{{\cal H}}$ is a
compact submanifold $M$ of $N_{{\cal H}}$ so that the inclusion of $M$
into $N_{{\cal H}}$ is a homotopy equivalence and so that $\partial M$
intersects every non-compact component of $\partial N_{{\cal H}}$ in
an incompressible annulus and contains every toroidal component of
$\partial N_{{\cal H}}$. (Results of McCullough \cite{mccullough} or
Kulkarni and Shalen \cite{kulkarni-shalen} guarantee that $N_{\cal H}$
has a relative compact core whenever $\pi_1(N)$ is finitely
generated.)  Set $P=\partial M\cap \partial 
N_{{\cal H}}$.  The {\em ends} of $N_{{\cal H}}$ are in one-to-one
correspondence with the components of $\partial M -P$ (see Bonahon
\cite{bonahon}).  An end $E$ of $N_{{\cal H}}$ is {\em geometrically
finite} if it has a neighborhood $U$ such that $U\cap C(N)
=\emptyset$, and is {\em geometrically infinite} otherwise. We also
refer to the corresponding components of $\partial M -P$ as
geometrically finite or geometrically infinite.

The following generalization of Thurston's covering theorem
\cite{thurston-notes} appears in \cite{canary-covering}. 

\begin{theorem}{covering-theorem}{}{} Let $N=\H^3/\Gamma$ be a
topologically tame hyperbolic $3$-manifold which covers another
hyperbolic 3-manifold $\hat N$ by a local isometry $\pi : N \to \hat
N$, and suppose that $\Omega(\Gamma)$ is non-empty.  If $E$ is a
geometrically infinite end of $N_{\cal H}$, then $E$ has a
neighborhood $U$ such that $\pi$ is finite-to-one on $U$.
\end{theorem}

The relative  compact core encodes much of the structure of the Kleinian
group. For example, a
finitely generated, torsion-free Kleinian group $\Gamma$ with
associated relative compact core $(M,P)$ has connected limit set if
and only if every geometrically finite component of $\partial M-P$ is
incompressible.  In this language, an accidental parabolic gives rise
to an essential annulus $A$ in $M$ which has one boundary component in
$P$ and the other in a geometrically finite component of $\partial
M-P$.  If $\Gamma$ has connected limit set, then the non-cyclic
subgroups in the Abikoff-Maskit decomposition of $\Gamma$ are simply
the fundamental groups of the components of $M-{\cal A}$, where ${\cal
A}$ is a maximal collection of disjoint, non-parallel essential annuli
associated to accidental parabolics.  If $\Gamma$ has disconnected
limit set, the non-cyclic subgroups in an Abikoff-Maskit decomposition of
$\Gamma$ into groups with connected limit set and elementary groups
arise as fundamental groups of components of $M-{\cal C}$,
where ${\cal C}$ is a maximal collection of disjoint, non-parallel
compressing disks for $(M,P)$ whose boundary curves lie in
geometrically finite components of $\partial M-P$. 

\section{The main results}
\label{strong limits}

In this section we prove our main result, namely
that a type-preserving sequence converges strongly, under the
assumption that the domain of discontinuity of the algebraic limit is
non-empty.

\begin{theorem}{strong-limit}{}{} Let $G$ be a finitely generated,
torsion-free group, and let $\{\rho_j\}\subset\calD(G)$ be a
type-preserving sequence converging algebraically to
$\rho\in\calD(G)$.  If $\Omega(\rho(G))$ is non-empty, then
$\{\rho_j\}$ converges strongly to $\rho$ and $\{\Lambda(\rho_j(G))\}$
converges to $\Lambda(\rho(G))$.
\end{theorem}

Before beginning the proof, we recall the statement of Theorem F from
\cite{anderson-canary}.

\begin{theorem}{}{}{(Theorem F from \cite{anderson-canary})}
Let $G$ be a finitely generated, torsion-free, non-abelian
group and let $\{\rho_j$\} be a sequence in ${\cal D}(G)$ converging to $\rho$.
If $\Lambda(\rho(G))=\rs$  and
$G$ is not a non-trivial free product of  (orientable)
surface groups and cyclic groups,
then $\{\rho_j\}$ converges strongly to $\rho$.
Moreover, $\{\Lambda(\rho_j(G))\}$ converges to $\Lambda(\rho(G))=\rs$.
\end{theorem}

Combining Theorem \ref{strong-limit} with the above theorem
yields the following immediate corollary.

\begin{corollary}{group-condition}{}{} Let $G$ be a finitely
generated, torsion-free group, and let $\{\rho_j\}\subset\calD(G)$ be
a type-preserving sequence converging algebraically to $\rho$.  If $G$
is not a (non-trivial) free product of (orientable) surface groups and
infinite cyclic groups, then $\{\rho_j\}$ converges strongly to $\rho$
and $\{\Lambda(\rho_j(G))\}$ converges to $\Lambda(\rho(G))$.
\end{corollary}

We further note that we obtain the same conclusion if we replace the
group-theoretic assumption in the statement of the corollary with the
weaker topological assumption that the  relative compact core of the
algebraic limit is not a relative compression body.  For a discussion
of the relation between the group-theoretic and topological
conditions, we refer the reader to Section $11$ of
\cite{anderson-canary}.

\begin{proof}{Theorem \ref{strong-limit}} The proof divides naturally
into three steps.  In Step $1$, we show that Theorem
\ref{strong-limit} holds in the case that $\rho(G)$ is either
a generalized web group, a degenerate group without accidental
parabolics or an elementary group.  In Step $2$, we show that Theorem
\ref{strong-limit} holds in the case that $\rho(G)$ has connected
limit set, using the result of Step 1 and the Abikoff-Maskit decomposition of
$\rho(G)$ into generalized web groups, degenerate groups without accidental
parabolic elements and elementary groups.  In Step $3$, we
pass from the case of connected limit set to the general case, 
making use of the Abikoff-Maskit decomposition of a torsion-free
finitely generated Kleinian group
into groups with connected limit sets and elementary groups.

\medskip
\noindent
{\bf Step 1:}  We recall that Thurston \cite{thurston-notes} and
Kerckhoff and Thurston \cite{kerckhoff-thurston} established Theorem
\ref{strong-limit} in the case that the algebraic limit is a
degenerate group without accidental parabolic elements.
One can obtain this by combining Theorem 9.6.1 in
Thurston \cite{thurston-notes} with the proof of Corollary 2.2 in
Kerckhoff and Thurston \cite{kerckhoff-thurston}; see also Corollary
6.1 in Ohshika \cite{ohshika-divergent}.

\medskip

We next establish Theorem \ref{strong-limit} for generalized web
groups.

\begin{proposition}{component-preserved}{}{} Let $G$ be a finitely
generated, torsion-free, non-abelian group, and let
$\{\rho_j\}\subset\calD(G)$ be a type-preserving sequence converging
algebraically to $\rho$. If $\rho(G)$ is a generalized web group, then
$\{\rho_j\}$ converges strongly to $\rho$ and $\{\Lambda(\rho_j(G))\}$
converges to $\Lambda(\rho(G))$.
\end{proposition}

\begin{proof}{Proposition \ref{component-preserved}} Let $\Gamma^0$ be
a geometrically finite subgroup of $\rho(G)$.  Since each
$\rho_j\circ\rho^{-1}$ is type-preserving, Marden's  quasiconformal
stability theorem (Proposition 9.1 in Marden \cite{marden}) implies
that $\Gamma^0_j=\rho_j\circ\rho^{-1}(\Gamma^0)$ is geometrically
finite and quasiconformally conjugate to $\rho(G)$ for sufficiently
large $j$. In particular, there exists a sequence $\{\phi_j\}$ of
$k_j$-quasiconformal maps of $\rs$ converging to the identity map,
with $\{ k_j\}$ converging to $1$, so that
$\rho_j\circ\rho^{-1}(\gamma) =\phi_j\circ\gamma\circ \phi_j^{-1}$ for
all $\gamma\in\Gamma^0$.  One sees immediately that $\{\Gamma_j^0\}$
converges geometrically to $\Gamma^0$ and that
$\{\Lambda(\Gamma^0_j)\}$ converges to $\Lambda(\Gamma^0)$. Hence,
Theorem \ref{strong-limit} holds whenever $\rho(G)$ is geometrically
finite, and in particular for $\rho(G)$ quasifuchsian or extended
quasifuchsian.

Suppose now that $\rho(G)$ is a web group.  Proposition \ref{limcon}
implies that it suffices to prove that $\{\Lambda(\rho_j(G)\}$ converges
to $\Lambda(\rho(G))$. If $\{\Lambda(\rho_j(G))\}$ does not converge to
$\Lambda(\rho(G))$, then there exists a subsequence of $\{\rho_j\}$, again
called $\{\rho_j\}$, so that $\{\Lambda(\rho_j(G))\}$ converges to a
set $\hat\Lambda$ which contains $\Lambda(\rho(G))$ as a proper
subset.  In particular, there must exist a quasifuchsian component
subgroup of $\rho(G)$ whose limit set separates $\hat\Lambda$. 
Proposition \ref{alg-geom} implies that we may pass to a further
subsequence of $\{\rho_j\}$, still called $\{\rho_j\}$, such that
$\{\rho_j(G)\}$ converges geometrically to a Kleinian group $\hat\Gamma$.
The definitions of geometric and Hausdorff convergence imply
that $\Lambda(\rho(G))\subset\Lambda(\hat\Gamma)\subset\hat\Lambda$.

Let $x\in \hat\Lambda-\Lambda(\rho(G))$, let $\Delta$ be the
component of $\Omega(\rho(G))$ which contains $x$, and set
$\Phi=\st_{\rho(G)}(\Delta)$, so that $\Lambda(\Phi)$ separates
$\hat\Lambda$.  Set $\Phi_j =\rho_j\circ\rho^{-1}(\Phi)$.  Since
$\rho(G)$ is a web group, $\Phi$ is quasifuchsian, and hence
geometrically finite. Let $\{ \phi_j\}$ be the sequence of
quasiconformal maps conjugating $\Phi_j$ to $\Phi$ produced in the
first paragraph of this proof.  Since $\Lambda(\Phi)$ separates
$\hat\Lambda$, $\{\Lambda(\rho_j(G))\}$ converges to $\hat \Lambda$,
and $\{\phi_j\}$ converges to the identity map, we see
that $\Lambda(\Phi_j)=\phi_j(\Lambda(\Phi))$
separates $\Lambda(\rho_j(G))$ for all sufficiently large $j$.

Recall that a {\em spanning disc} for the quasifuchsian subgroup
$\Phi$ of the Kleinian group $\rho(G)$ is a properly embedded open
disc $D$ in $\bfH^3$ which is precisely invariant under $\Phi$ in
$\rho(G)$ and which extends to a closed disc $\overline{D}$ in
$\bfH^3\cup\rs$ with boundary $\Lambda(\Phi)$.
Let $N =\bfH^3/\rho(G)$ and $\hat N =\bfH^3/\hat\Gamma$, and
$\pi: N\rightarrow\hat N$ and $\hat p: \bfH^3\rightarrow\hat N$ be the
associated covering maps.  It is shown in the
proof of Proposition 6.1 of \cite{anderson-canary} that, for any
algebraically convergent sequence in $\calD(G)$ whose algebraic limit
is a web group, such as $\rho(G)$, and for any component subgroup of
$\rho(G)$, such as $\Phi$, there exists a spanning disk $D$ for $\Phi$
and a compact core $M$ for $N$ such that $\pi$ is an embedding
restricted to $M$, $\hat p(D)$ is a properly embedded surface in $\hat
N$, and $\pi(M)$ is disjoint from $\hat p(D)$.

Set $N_j =\bfH^3/\rho_j(G)$.  Let $\{ V_j\subset N_j\}$ and $\{
f_j:V_j\to \hat N\}$ be the sequence of submanifolds and
$K_j$-bilipschitz diffeomorphisms produced by Lemma \ref{fi
convergence}.  Since $\pi(M)$ lies in $f_j(V_j)$ for all sufficiently
large $j$, the submanifold $M_j =f_j^{-1}(\pi(M))$ of $N_j$ is
defined.  In Corollary C of \cite{anderson-canary}, we observe that
$M_j$ is a compact core for $N_j$ for all sufficiently large $j$.

One may extend each $\phi_j$ to a $L_j$-bilipschitz diffeomorphism
$\psi_j:\bfH^3\to\bfH^3$ conjugating the action of $\Phi$ to the
action of $\Phi_j$ (see Douady and Earle \cite{douady-earle}, Reimann
\cite{reimann}, or Tukia \cite{tukia}) in such a way that $\{ L_j\}$
converges to $1$ and $\{\psi_j\}$ converges to the identity map.  Let
$D_j=\psi_j(D)$, and note that $D_j$ is a properly embedded open disc
in $\bfH^3$ which is invariant under $\Phi_j$ and which extends to a
closed disc in $\bfH^3\cup\rs$ with boundary $\Lambda(\Phi_j)$.
However, $D_j$ need not be a spanning disk for $\Phi_j$, as it need
not be precisely invariant under $\Phi_j$ in $\rho_j(G)$.

We next show that $p_j(D_j)$ is disjoint from $M_j$ for sufficiently
large $j$, where $p_j:\bfH^3\rightarrow N_j$ is the covering map.
Choose $\delta>0$ so that $\min_{x\in\pi(M)}\inj_{\hat N}(x)>3\delta$.
Since $f_j: V_j\rightarrow\hat N$ is $K_j$-bilipschitz and $\{ K_j\}$
converges to $1$, we see that $\min_{x\in M_j}\inj_{N_j}(x) > 2\delta$
for sufficiently large $j$.

Let $X\subset D$ denote the set of points for which there exists a
non-trivial $\gamma\in\Phi$ such that $\rmd(x,\gamma(x))< \delta$, and
let $Y=D-X$.  Consider the subsets $Y_j=\psi_j(Y)$ and $X_j=\psi_j(X)$
of $D_j$.  For each $x\in X_j$, there exists a non-trivial
$\gamma_j\in\Phi_j$ so that $\rmd(x, \gamma_j(x)) < L_j\delta$.
Since  $L_j\delta < \min_{x\in M_j}\inj_{N_j}(x)$ for
sufficiently large $j$, we have that $p_j(X_j)$ is disjoint from
$M_j$.  Moreover, since $\hat p(Y)$ is compact and is disjoint from
$\pi(M)$, and since both $\{ f_j\}$ and $\{\psi_j\}$ converge to the
identity map, $p_j(Y_j)$ is disjoint from $M_j$ for sufficiently large
$j$.  This completes the proof that $p_j(D_j)$ is disjoint from $M_j$
for sufficiently large $j$.

We now assume that we have chosen $j$ large enough so that
$\Lambda(\Phi_j)$ separates $\Lambda(\rho_j(G))$,  $M_j$ is a
compact core for $N_j$, and $p_j(D_j)$ is disjoint from $M_j$.
Since $\Lambda(\Phi_j)$ separates $\Lambda(\rho_j(G))$ and since pairs of
fixed points of loxodromic elements of $\rho_j(G)$ are dense in
$\Lambda(\rho_j(G))\times\Lambda(\rho_j(G))$ (see Proposition
V.E.5 in \cite{maskit-book}), there exists a
hyperbolic element $\gamma_j\in \rho_j(G)$ whose fixed points 
are separated by $\Lambda(\Phi_j)$.  If $C_j$ is the closed geodesic
in $N_j$ which is the projection of the axis of $\gamma_j$, then every
curve homotopic to $C_j$ must intersect $p_j(D_j)$. However, since
$M_j$ is a compact core for $N_j$ which is disjoint from $p_j(D_j)$,
there exists a curve in $M_j$ which is homotopic to $C_j$ and disjoint
from $p_j(D_j)$.  This contradiction establishes the result. 
\end{proof}

We note the following corollary of the proof of Proposition
\ref{component-preserved}.

\begin{corollary}{generalized-web}{}{} Let $G$ be a finitely
generated, torsion-free, non-abelian group, and let
$\{\rho_j\}\subset\calD(G)$ be a type-preserving sequence converging
algebraically to $\rho$.  Suppose that $\rho(G)$ is a generalized web
group, and let $\Phi$ be a component subgroup of $\rho(G)$.  Then,
$\Phi_j =\rho_j\circ\rho^{-1}(\Phi)$ is quasifuchsian and is a
component subgroup of $\rho_j(G)$ for all sufficiently large $j$.
\end{corollary}

We complete Step 1 by establishing Theorem \ref{strong-limit} for
elementary groups, under the additional assumption that the limit
representation $\rho$ is discrete and faithful.

\begin{lemma}{elementary}{}{} Let $G$ be a free abelian group of rank
at most two, and let $\{\rho_j\}$ be a type-preserving sequence in
$\calD(G)$ converging algebraically to $\rho\in\calD(G)$.  Then,
$\{\rho_j\}$ converges strongly to $\rho$.
\end{lemma}

\begin{proof}{Lemma \ref{elementary}} Suppose there exists a sequence
of non-trivial elements $\{ g_j\}$ of $G$ so that $\{\rho_j(g_j)\}$
converges to a M\"obius transformation $\gamma$.  We wish to show that
$\gamma\in\rho(G)$.  We begin by noting that $\fix(\gamma)
=\Lambda(\rho(G))$, for otherwise, the fixed point set
$\fix(\rho_j(g_j))$ of $\rho_j(g_j)$ would be distinct from
$\Lambda(\rho_j(G))$ for all sufficiently large $j$, which cannot
occur.

Suppose that $G =\langle g\rangle$ is infinite cyclic and that
$\rho(g)$ is loxodromic.  As $G$ is cyclic, we may write $g_j
=g^{n_j}$ for $n_j\in\bfZ$. Since $\{\rho_j(g)\}$ converges to
$\rho(g)$, and since $\rho_j(g)$ is loxodromic for all $j$ and
$\rho(g)$ is loxodromic, there exists a sequence of M\"obius
transformations $\{\beta_j\}$ converging to the identity so that the
axis of $\beta_j\rho_j(g)\beta_j^{-1}$ in $\bfH^3$ is equal to the
axis $l$ of $\rho(g)$ in $\bfH^3$ for all $j$.  Since $\{\beta_j\}$
converges to the identity, $\{\beta_j\rho_j(g_j)\beta_j^{-1}\}$
converges to $\gamma$, and so the translation length of
$\beta_j\rho_j(g_j)\beta_j^{-1}$ along $l$ converges to the
translation length of $\gamma$.

If $\{ n_j\}$ is not bounded, there exists a subsequence of $\{
n_j\}$, again called $\{ n_j\}$, converging to either $\infty$ or
$-\infty$.  Since the translation length of $\rho_j(g_j)$ is $n_j$
times the translation length of $\beta_j\rho_j(g)\beta_j^{-1}$, we see
that the translation length of $\beta_j\rho_j(g)\beta_j^{-1}$ along
$l$ converges to zero, and so $\rho(g)$ fixes $l$ pointwise.  In
particular, $\rho(g)$ is either the identity or elliptic.  However,
$\rho(g)$ cannot be the identity, as this would contradict the
faithfulness of $\rho$.  Also, $\rho(g)$ cannot be elliptic, as this
would contradict either the faithfulness of $\rho$ if $\rho(g)$ has
finite order or the discreteness of $\rho(G)$ if $\rho(g)$ has
infinite order. This contradiction gives that $\{ n_j\}$ is bounded.
If $\{ n_j\}$ is not eventually constant, there exist two subsequences
with limits $m\neq n$, and so there are two subsequences of
$\{\rho_j(g_j)\}$ whose limits are $\rho(g)^m$ and $\rho(g)^n$, which
cannot occur. Hence, $\{ n_j\}$ is eventually constant, and so
$\gamma\in\rho(G)$. 

The proof in the case that $\rho(G)$ is a parabolic subgroup is
similar.  Suppose that $G$ is infinite cyclic and $\rho(g)$ is
parabolic.  Since $\{\rho_j(g)\}$ converges to $\rho(g)$, and since
$\rho_j(g)$ is parabolic for all $j$ and $\rho(g)$ is parabolic, there
exists a sequence of M\"obius transformations $\{\beta_j\}$ converging
to the identity so that $\beta_j\rho_j(g)\beta_j^{-1} =\rho(g)$ for
all $j$.  Normalizing so that $\rho(g)(z)=z+1$ and writing $g_j
=g^{n_j}$ for $n_j\in\bfZ$, we see that 
$\beta_j\rho_j(g_j)\beta_j^{-1}(z)=z+n_j$.  Since
$\{\beta_j\rho_j(g_j)\beta_j^{-1}\}$ converges to a M\"obius
transformation $\gamma$, 
it is easy to see that $\{ n_j\}$ must be eventually constant, and
so $\gamma\in\rho(G)$. 

In the case that $G$ has rank two, choose generators $a$ and $b$.
Since $\rho_j(a)$ and $\rho_j(b)$ are parabolic for all $j$, and since
$\rho(a)$ and $\rho(b)$ are parabolic, there exists a sequence of
M\"obius transformations $\{\beta_j\}$ converging to the identity so
that $\beta_j\rho_j(a)\beta_j^{-1} =\rho(a)$ for all $j$.  Normalizing
so that $\rho(a)(z) =z+1$, we see that
$\beta_j\rho_j(b)\beta_j^{-1}(z)=z+\tau_j$ and $\rho(b)(z)=z+\tau$,
where $\{\tau_j\}$ converges to $\tau$ and $\tau$ has non-zero
imaginary part.  
Write $g_j =a^{n_j} b^{m_j}$ for $n_j$, $m_j\in\bfZ$, so that
$\beta_j\rho_j(g_j)\beta_j^{-1}(z)=z+n_j+m_j\tau_j$.  Since
$\{\beta_j\rho_j(g_j)\beta_j^{-1}\}$ converges to a M\"obius
transformation $\gamma$, it is easy
to see that both $\{ n_j\}$ and $\{ m_j\}$ must be eventually
constant, and so $\gamma\in\rho(G)$. 
\end{proof}

\medskip
\noindent
{\bf Step 2:}  In this step, we establish Theorem \ref{strong-limit}
in the case in which the algebraic limit $\rho(G)$ has connected limit
set.  As we have already shown that Theorem 3.1 holds for elementary groups,
generalized web groups,  and degenerate groups without accidental parabolics,
it  will suffice to show that it is
preserved by Klein-Maskit combination along a parabolic cyclic subgroup.

In order to cleanly handle comparisons of cusp regions of the limit of
a convergent sequence and of the approximates, we use the fact that we
are working with type-preserving sequences to make a convenient
normalization.  Let $\{\rho_j\}\subset\calD(G)$ be a type-preserving
sequence converging algebraically to $\rho$, and let $\rho(h)$ be a
primitive parabolic element of $\rho(G)$. Since $\{\rho_j(h)\}$ 
is a sequence of parabolic M\"obius transformations converging to
$\rho(h)$, there exists a sequence $\{\beta_j\}$ of M\"obius
transformations converging to the identity so that 
$\beta_j\rho_j(h)\beta_j^{-1}(z)=\rho(h)$ for all $j$.  Since
$\{\beta_j\}$ converges to the identity, the modified sequence
$\{\beta_j\rho_j\beta_j^{-1}\}$ has the same algebraic limit as the
original sequence $\{\rho_j\}$, namely $\rho$, and the geometric limit
of $\{\beta_j\rho_j(G)\beta_j^{-1}\}$ is equal to the geometric limit
of the original sequence $\{\rho_j(G)\}$, if the geometric limit
exists.  Hence, we may replace the original sequence with the modified
sequence without effect. We refer to this process as {\em normalizing
the sequence about $\rho(h)$}.

We now wish to describe what we mean by cusp regions persisting in the
approximates.  Let $\{\rho_j\}\subset\calD(G)$ be a type-preserving
sequence converging to $\rho$ and let $D\subset\Omega(\rho(G))$ be a
cusp region for $\rho(h)$.  Say that $D$ {\em persists in the
approximates} if, whenever $\{\rho_j\}$ is normalized about $\rho(h)$
to obtain a new sequence $\{\rho_j'\}$, $D$ is a cusp region for
$\rho_j'(h)$ for all sufficiently large $j$.

In order to show that the conclusions of Theorem \ref{strong-limit}
are preserved under the application of the Klein-Maskit combination
theorems, we need to carry along the information about the persistence
of cusp regions.  We begin by showing in Proposition \ref{persistence}
that cusp regions persist in the approximates for degenerate groups
without accidental parabolic elements and for generalized web groups.
In Proposition \ref{accidental}, we show that if $\rho(G)$ is obtained
from $\rho(G_1)$ and $\rho(G_2)$ by Klein-Maskit combination along an
accidental parabolic element, and if $\{\Lambda(\rho_j(G_m)\}$
converges to $\Lambda(\rho(G_m)$ and all cusp regions of $\rho(G_m)$
persist in the approximates for both $m$, then
$\{\Lambda(\rho_j(G))\}$ converges to $\Lambda(\rho(G))$, $\{\rho_j\}$
converges strongly to $\rho$, and all cusp regions of $\rho(G)$
persist in the approximates.  Combining these two propositions, we
construct an inductive proof of Theorem \ref{strong-limit} in the case
that $\rho(G)$ has connected limit set.

\begin{proposition}{persistence}{}{} Let $G$ be a finitely generated,
torsion-free group, and let $\{\rho_j\}\subset\calD(G)$
be a type-preserving sequence converging algebraically to $\rho\in {\cal D}(G)$.
Suppose that $\rho(G)$ is either degenerate without
accidental parabolics, generalized web, or elementary.  Then, cusp
regions persist in the approximates.
\end{proposition}

\begin{proof}{Proposition \ref{persistence}} We begin by noting that
we have already shown that $\{\rho_j\}$ converges strongly to $\rho$
and that $\{\Lambda(\rho_j(G))\}$ converges to $\Lambda(\rho(G))$.
Let $D$ be a cusp region for the primitive parabolic $\rho(h)$. We
assume that the sequence $\{\rho_j\}$ is normalized about $\rho(h)$.
In the case that $\rho(G)$ is elementary, there is nothing to prove.
Hence, we may assume that $\rho(G)$ is either degenerate without
accidental parabolics or generalized web.

Let $c =\partial D \cap \Omega(\rho(G))$.  The first key observation
is that $c$ lies in $\Omega(\rho_j(G_m))$ and is precisely invariant
under $\langle\rho_j(h)\rangle$ in $\rho_j(G)$ for sufficiently large
$j$.  As this argument is used later, we state it as a lemma.

\begin{lemma}{picusps}{}{} Let $G$ be a finitely generated,
torsion-free, non-abelian group, and let $\{\rho_j\}\subset\calD(G)$
be a type-preserving sequence converging strongly to $\rho$ such that
$\{\Lambda(\rho_j(G))\}$ converges to $\Lambda(\rho(G))$.  Suppose
that $D$ is a cusp region for $\rho(G)$ associated to the primitive
parabolic element $\rho(h)$, that $\{\rho_j\}$ is normalized about
$\rho(h)$, and set $c=\partial D\cap \Omega(\rho(G))$.  Then, for 
sufficiently large $j$, $c$ lies in $\Omega(\rho_j(G))$ and is
precisely invariant under $\langle\rho_j(h)\rangle$ in $\rho_j(G)$.
\end{lemma}

\begin{proof}{} We may assume that $\rho(h)(z)=z+1$.  Choose a point
$z$ on $c$, and let $c_f$ be the closed arc in $c$ joining $z$ to
$z+1$, so that the interior of $c_f$ is a fundamental domain for the
action of $\langle\rho(h)\rangle$ on $c$.  Since $c_f$ is compact and
since $\{\Lambda(\rho_j(G))\}$ converges to $\Lambda(\rho(G))$, we
see that $c_f$ lies in $\Omega(\rho_j(G))$ for sufficiently large
$j$, and hence that $c$ lies in $\Omega(\rho_j(G))$ for sufficiently
large $j$.

In order to show that $c$ is precisely invariant under
$\langle\rho_j(h)\rangle$ in $\rho_j(G)$, we need to make use of a
standard {\em convergence property} for M\"obius transformations,
which we describe here.  Let $K$ be a compact set in $\rs$ and let
$\{ M_j\}$ be a sequence of distinct M\"obius transformations so
that $K\cap M_j(K)$ is non-empty for all $j$.  Then, there exists a 
subsequence of $\{ M_j\}$, again called $\{ M_j\}$, so that either
there exists a sequence of points $\{ x_j\in\fix(M_j)\}$ converging to
a point in $K$ or $\{ M_j\}$ converges to a M\"obius transformation
(see Marden \cite{marden-survey} or Gehring and Martin
\cite{gehring-martin}.)

If $c$ is not precisely
invariant under $\langle\rho_j(h)\rangle$ in $\rho_j(G)$ for all
sufficiently large $j$, there exists a sequence $\{ g_j\}$ of elements of
$G -\langle h\rangle$ so that $c\cap\rho_j(g_j)(c)$ is non-empty for infinitely
many $j$.  We
may pre- and post-multiply each $g_j$ by appropriate powers of $h$ to
produce a sequence $\{ g_j'\}$ of elements of $G -\langle h\rangle$ so
that $c_f\cap\rho_j(g_j')(c_f)$ is non-empty for infinitely many $j$. The
convergence property for M\"obius transformations implies that there
exists a subsequence of $\{\rho_j(g'_j)\}$, again called
$\{\rho_j(g'_j)\}$, for which  either there exist points $\{ x_j\in
\Lambda(\rho_j(G))\}$ converging to a point in $c_f$, or
$\{\rho_j(g_j')\}$ converges to a M\"obius transformation $\gamma$. In
the latter case, note that $c_f\cap\gamma(c_f)$ is necessarily 
non-empty. 

The former case cannot occur, as $\{\Lambda(\rho_j(G))\}$ converges
to $\Lambda(\rho(G))$, which is disjoint from $c_f$.  The latter
case cannot occur, as the precise invariance of $c$ under
$\langle\rho(h)\rangle$ in $\rho(G)$ and the strong convergence of $\{
\rho_j\}$ to $\rho$ together imply that $\gamma =\rho(h^n)$ for some
$n\ne 0$; Lemma 3.6 of J\o rgensen and Marden \cite{jorgensen-marden}
then implies that $g_j'=h^n$ for all sufficiently large $j$, contradicting
our assumption that $g_j'$ does not lie in $\langle h\rangle$.
\end{proof}

We now show that $D$ persists in the approximates.  Suppose first that
$\rho(G)$ is a degenerate group without accidental parabolics.
As each isomorphism
$\rho_j\circ\rho^{-1}:\rho(G)\rightarrow\rho_j(G)$ is type-preserving and
$\Omega(\rho_j(G))$ is non-empty for all large enough $j$,
Theorem 6 of Maskit \cite{maskit-koebe} gives that,
for all large enough $j$,
$\rho_j(G)$ does not contain accidental parabolic elements.
However, since $c$ is precisely invariant under $\langle\rho_j(h)\rangle$ in
$\rho_j(G)$ for all sufficiently large $j$, by Lemma \ref{picusps}, we see
that either $D$ is a cusp region for
$\rho_j(h)$ in $\rho_j(G)$, or $C =c\cup {\rm fix}(\rho(h))$
 separates $\Lambda(\rho_j(G))$, in
which case $\rho_j(h)$ is an accidental parabolic element of
$\rho_j(G)$.  As the latter case cannot occur, it must be that $D$ is
a cusp region for $\rho_j(h)$ in $\rho_j(G)$ for all sufficiently large $j$.


Suppose now that $\rho(G)$ is a generalized web group, let $\Delta$ be
the component of $\Omega(\rho(G))$ containing $D$, and let $\Phi
=\st_{\rho(G)}(\Delta)$.  By Lemma \ref{picusps}, we know that $c$
lies in $\Omega(\rho_j(G))$ and is precisely invariant under
$\langle\rho_j(h)\rangle$ for all sufficiently large $j$.  By
Corollary \ref{generalized-web}, $\Phi_j =\rho_j\circ\rho^{-1}(\Phi)$
is a quasifuchsian component subgroup of $\rho_j(G)$, stabilizing the
component $\Delta_j$ of $\rho_j(G)$.  Since $\{\Lambda(\Phi_j)
=\partial\Delta_j\}$ converges to $\Lambda(\Phi) =\partial\Delta$, we
see that $c$ is contained in $\Delta_j$ for $j$ sufficiently
large.  Since quasifuchsian component subgroups cannot contain
accidental parabolic elements (see Theorem 4 of Maskit
\cite{maskit-boundaries}), the Jordan curve $C =\partial D$ cannot
separate  $\Lambda(\rho_j(G))$, and so $D$ is a cusp region for
$\rho_j(h)$ in $\rho_j(G)$ for all sufficiently large $j$.
\end{proof}

The following proposition shows that the conclusions of Theorem
\ref{strong-limit} are preserved by Klein-Maskit combination along
a cyclic parabolic subgroup.

\begin{proposition}{accidental}{}{} Let $G$ be a finitely generated,
torsion-free, non-abelian group, and let $\{\rho_j\}\subset\calD(G)$
be a type-preserving sequence converging algebraically to $\rho$ with
$\Omega(\rho(G))$ non-empty.  Suppose that $\Theta
=\langle\rho(h)\rangle$ is an accidental parabolic subgroup of
$\rho(G)$, and that $\rho(G_1)$ and $\rho(G_2)$ are the factor
subgroups of the Klein-Maskit decomposition of $\rho(G)$ along
$\Theta$.  If $\{\Lambda(\rho_j(G_m))\}$ converges to
$\Lambda(\rho(G_m))$ for both $m$ and if all cusp regions for
$\rho(G_m)$ persist in the approximates, then $\{\rho_j\}$ converges
strongly to $\rho$, $\{\Lambda(\rho_j(G))\}$ converges to
$\Lambda(\rho(G))$, and all cusp regions for $\rho(G)$ persist in the
approximates.
\end{proposition}

\begin{proof}{Proposition \ref{accidental}}  We begin by showing that 
$\{\Lambda(\rho_j(G))\}$ converges to $\Lambda(\rho(G))$.
Proposition \ref{limcon} then implies the strong convergence of
$\{\rho_j\}$ to $\rho$.  In order to show that the limit sets
converge, we first show that the approximates are also
Klein-Maskit combinations along (the same) cyclic parabolic subgroup.
We use this to construct explicit fundamental domains for the approximates
which can be used to show that the limit sets converge.
We then show that the cusp regions for
$\rho(G)$ persist in the approximates.  The details for the two
Klein-Maskit combinations are similar; we give full details for the
case of a Klein-Maskit combination of type I, and briefly sketch the
argument for a Klein-Maskit combination of type II at the end of the
proof.

If $\{\Lambda(\rho_j)\}$ does not converge to $\Lambda(\rho(G))$, we
may pass to a subsequence so that $\{\Lambda(\rho_j(G))\}$ converges to a
closed set $\hat\Lambda$ containing $\Lambda(\rho(G))$ as a proper subset.
Without loss
of generality, we may assume that $\{\rho_j\}$ is normalized about
$\rho(h)$ and that $\rho(h)(z)=z+1$.  Let $c$ be the axis of $\rho(h)$
and let $C = c\cup\{\infty\}$ be the full axis of $\rho(h)$.  Let
$E_1$ and $E_2$ be the two closed discs in $\rs$ determined by $C$,
labeled so that $E_m$ is a cusp region for $\rho(h)$ in $\rho(G_m)$. 

Since cusp regions persist in the approximates for both $\rho(G_m)$,
we have that $E_m$ is a cusp region for $\rho_j(h)$ in $\rho_j(G_m)$
for both $m$ and for all sufficiently large $j$.
Therefore, $\rho_j(G)$ is formed from $\rho_j(G_1)$ and $\rho_j(G_2)$ by a
Klein-Maskit combination of type I along $C$ for all sufficiently
large $j$. 

\medskip

Choose a point $z$ on $c$ and
consider the open vertical band $V$ in $\bfC$ whose two bounding lines
pass through $z$ and $z+1$.  For both $m$, we can choose a fundamental
domain $F^m$ for $\rho(G_m)$ which lies in $V$ and contains $V\cap
E_m$.  Since $\rho(G)$ is formed from $\rho(G_1)$ and $\rho(G_2)$ by
Klein-Maskit combination of type I along $C$, Theorem \ref{kmcombI}
assures us that $F =F^1\cap F^2$ is a fundamental domain for
$\rho(G)$.  Choose a point $x\in \hat\Lambda-\Lambda(\rho(G))$ which
lies in the interior of $F$ (noting that one may need to alter the
choice of $F$ to guarantee that $x$ lies in $F$ rather than on its
boundary).  Let $K$ be a compact set which contains an open
neighborhood of $x$ and which is contained entirely in the interior
of $F$.  Since $K$ lies in $F$, it is precisely invariant under the
identity in $\rho(G)$.

By construction, $K$ lies in the strip $V$.  Combining the convergence
property for M\"obius transformations (see the proof of Lemma
\ref{picusps}) with the assumptions that $\{\rho_j(G_m)\}$ converges
to $\rho(G_m)$ and $\{\Lambda(\rho_j(G_m))\}$ converges to
$\Lambda(\rho(G_m))$, we see that $K$ lies in $\Omega(\rho_j(G_m))$
and is precisely invariant under the identity in $\rho_j(G_m)$ for
sufficiently large $j$.  Moreover, since $E_m$ is precisely invariant
under $\langle\rho_j(h)\rangle$ in $\rho_j(G_m)$ for sufficiently
large $j$, we may construct a fundamental domain $F^m_j$ for
$\rho_j(G_m)$ containing both $V\cap E_m$ and $K$ for sufficiently
large $j$; in particular, $F^1_j\cap c =F^2_j\cap c$.  Since
$\rho_j(G)$ is formed from $\rho_j(G_1)$ and $\rho_j(G_2)$ by
Klein-Maskit combination of type I along $C$, Theorem \ref{kmcombI}
implies that $F_j=F_j^1\cap F_j^2$ is then a fundamental domain for
$\rho_j(G)$ for sufficiently large $j$. 

Since $\{\Lambda(\rho_j(G))\}$ converges to $\hat\Lambda$, there
exists a sequence of points $\{ x_j\in \Lambda(\rho_j(G))\}$
converging to $x$.  This implies that $x_j\in K\subset
F_j\subset\Omega(\rho_j(G))$ for sufficiently large $j$, a
contradiction.  Hence, $\{\Lambda(\rho_j(G))\}$ converges to
$\Lambda(\rho(G))$, and so $\{\rho_j\}$ converges strongly to $\rho$.

\medskip

In order to complete the proof in the case of a Klein-Maskit
combination of type I, it remains only to argue that cusp regions for
$\rho(G)$ persist in the approximates.  Let $D'$ be a cusp region
for a primitive parabolic element $\rho(g)$ of $\rho(G)$.
After suitably normalizing, we first use Klein-Maskit Combination
Theorem \ref{kmcombI} and our assumption that cusps persist
in the approximates of the factors to show that a cusp region $D''$
contained in $D'$ (and missing all translates of one of the
cusp regions we are combining along)
does persist in the approximates. We then use the convergence property
and the fact that limit sets converge to show that $D'$ itself must
persist in the approximates.

Let $\{\beta_j\}$ be a
sequence of M\"obius transformations converging to the identity so
that $\beta_j\rho_j(g)\beta_j^{-1}=\rho(g)$ for all $j$.  Set 
$\rho_j'=\beta_j\circ\rho_j\circ\beta_j^{-1}$, and let
$c_j=\beta_j(c)$ and $E^j_m=\beta_j(E_m)$ for all $j$ and for both
$m$.  Notice that, for all large enough $j$,
$E^j_m$ is a cusp region for $\rho_j'(h)$ in
$\rho_j'(G_m)$ and $\rho_j'(G)$ is a Klein-Maskit combination of
type I of $\rho_j'(G_1)$ and $\rho_j'(G_2)$ along
$\langle\rho_j'(h)\rangle$ with associated axis $c_j$.  Moreover,
$\{\rho_j'\}$ converges strongly to $\rho$ and
$\{\Lambda(\rho_j'(G))\}$ converges to $\Lambda(\rho(G))$.

Let $\xi$ be the fixed point of $\rho(g)$, and let $d'=\partial D'-\{ \xi\}$.
Lemma \ref{picusps} implies that $d'$ is precisely invariant
under $\langle\rho_j(g)\rangle$ in $\rho_j(G)$ for sufficiently large
$j$.  By Theorem \ref{kmcombI}, every cusp region for $\rho(G)$ is
a cusp region for either $\rho(G_1)$ or $\rho(G_2)$, and
so we may assume that $D'$ is a cusp region for $\rho(G_1)$.  Since we
have assumed that cusp regions for $\rho_j(G_1)$ persist in the
approximates, we see that $D'$ is a cusp region for $\rho_j'(G_1)$ for
all sufficiently large $j$.  Moreover, since $E_1$ is not a cusp region
for $\rho(G)$, we see that $D'$ and $E_1$ are inequivalent cusp
regions for $\rho(G_1)$, and hence that $D'$ and $E_j^1$ are
inequivalent cusp regions for $\rho_j(G_1)$ for all sufficiently large $j$.

Let $D''\subset D'$ be a cusp region for $\rho(g)$
in $\rho(G_1)$ which does not intersect any $\rho(G_1)$-translate of
$E_1$, and suppose that there exists a $\rho_j'(G_1)$-translate of
$E_1^j$ which intersects $D''$ for infinitely many $j$. Since $D''$
and $E_1^j$ are inequivalent for all sufficiently large $j$, $D''$ cannot
contain or be contained in any translate of $E_1^j$.  Hence,
$d''=\partial D''\cap\Omega(\rho_j'(G_1))$ must intersect some
$\rho_j'(G_1)$-translate of
$e''=\partial E_1^j\cap\Omega(\rho_j'(G_1))$ for infinitely many $j$.  This,
however, can be ruled out by the argument in the proof of Lemma
\ref{picusps}.  Hence, for all sufficiently large $j$, $D''$ does not
intersect any $\rho_j'(G_1)$-translate of $E^j_1$, and so Theorem
\ref{kmcombI} guarantees that $D''$ is a cusp region for $\rho_j'(G)$
for sufficiently large $j$.  That is, $D''$ persists in the
approximates of $\rho(G)$.

Suppose that $D'$ itself does not persist in the approximates of $\rho(G)$.
Then, for infinitely many $j$, there exists
$g_j\in G-\langle g\rangle$ such that $\rho_j'(g_j)(D')$ intersects
$D'$.  Since $d'$ is precisely invariant under
$\langle\rho_j(g)\rangle$ in $\rho_j(G)$ for sufficiently large $j$,
this implies that $\rho_j(g_j)(\xi)$ is contained in the interior of
$D'$ for infinitely many $j$.  However, the argument in the proof of
Lemma \ref{picusps} shows that $D'-D''\subset \Omega(\rho_j(G))$ for all
sufficiently large $j$. (The key point is that there is a
compact fundamental domain
for the action of $\langle\rho(g)\rangle$ on $D'-int(D'')$, where $int(D'')$
is the interior of $D''$.)  Since $D''$ is a cusp region for $\rho_j(G)$
for all sufficiently large $j$, we see that the interior of $D'$ must be
contained in $\Omega(\rho_j(G))$ for all sufficiently large $j$, a
contradiction.  Hence, $D'$ is a cusp region for $\rho_j'(G)$ for
all sufficiently large $j$, and we have completed the proof of Proposition
\ref{accidental} in the case that $\rho(G)$ is formed from $\rho(G_1)$
and $\rho(G_2)$ by a Klein-Maskit combination of type I.

\medskip

The proof in the case of Klein-Maskit combination of type II is quite
similar.  Suppose that the factor subgroups of the decomposition are
$\rho(H)$ and $\langle\rho(g)\rangle$.  Let $C_1$, $C_2$, $E_1$,
$E_2$, $\theta_1=\rho(h_1)$ and $\theta_2=\rho(h_2)$ be as in the
statement of the Klein-Maskit combination theorem of type II. We may
assume that our sequence is normalized about $\rho(h_1)$. Let
$c_m=C_m\cap \Omega(H)$. Let $c_1^j=c_1$, $E_1^j=E_1$,
$c_2^j=\rho_j(g)(c_1)$ and $E_2^j=\rho_j(g)(E_2)$.  By assumption, for
sufficiently large $j$, $E_m^j$ is a cusp region for $\rho_j(h_m)$ in
$\rho_j(H)$.  As before, we see that $\rho_j(G)$ is a Klein-Maskit
combination of type II with factor subgroups $\rho_j(H)$ and
$\langle\rho_j(g)\rangle$, for sufficiently large $j$, and that
$\{\Lambda(\rho_j(G))\}$ converges to $\Lambda(\rho(G))$.  This again
implies that $\{\rho_j\}$ converges strongly to $\rho$.  The proof
that cusp regions persist in the approximates is also much as above. 
\end{proof}

We are now prepared to establish Theorem \ref{strong-limit} in the
case when $\rho(G)$ has connected limit set.  We suppose that
$\{\rho_j\}\subset {\cal D}(G)$ is a type-preserving sequence
converging algebraically to $\rho$, that $\Omega(\rho(G))$ is
non-empty, and that $\Lambda(\rho(G))$ is connected. Let
$\{\Gamma_1,\ldots, \Gamma_k,\Gamma_{k+1},\ldots, \Gamma_l\}$ be the
Abikoff-Maskit decomposition of $\rho(G)$ as described in Section
\ref{decomp}, so that each $\Gamma_i$ is either 
degenerate without accidental parabolics, generalized web or cyclic.
Let $H_i =\rho^{-1}(\Gamma_i)$.
By the results of Step 1, each of the restricted sequences 
$\{\rho_j|_{H_i}\}$ converges strongly to $\rho|_{H_i}$ and
$\{\Lambda(\rho_j(H_i))\}$ converges to $\Lambda(\rho(H_i))$.
Proposition \ref{persistence} implies that all cusp regions of
$\rho(H_i)$ persist in the approximates.  We set $\Gamma^1=\Gamma_1$,
and recall that $\Gamma^i$ is formed from $\Gamma^{i-1}$ and
$\Gamma_i$ by Klein-Maskit combination along a parabolic cyclic
subgroup, of type I for $1\le i\le k$ and of type II for $k<i$, and
$\rho(G) =\Gamma^l$. Let $H^i=\rho^{-1}(\Gamma^i)$.
Applying Proposition \ref{accidental}
inductively shows, for all $i$, that $\{\rho_j|_{H^i}\}$ converges
strongly to $\rho|_{H^i}$, that $\{\Lambda(\rho_j(H^i))\}$ converges
to $\Lambda(\rho(H^i))$, and that all cusp regions of $\rho(H^i)$
persist in the approximates.  This completes the proof in the case
that $\rho(G)$ has connected limit set.

\medskip
\noindent
{\bf Step 3:} We conclude the proof of the theorem by considering the
case in which $\rho(G)$ has non-connected limit set.  In this case we
may decompose $\rho(G)$ into a finite collection of subgroups,
each of which either has connected limit set or is elementary, so
that $\rho(G)$ is built up from the subgroups by repeatedly applying
the operation of Klein combination. One may then complete the proof
of the theorem by applying the following result from 
\cite{anderson-canary} finitely many times.

\begin{proposition}{klein-convergence}{}{(Proposition 10.2
in \cite{anderson-canary})}{}
Let $G$ be a finitely generated, torsion-free group, and let
$\{\rho_j\}\subset\calD(G)$ be a type-preserving sequence converging
algebraically to $\rho\in\calD(G)$. Suppose that $\Omega(\rho(G))$ 
is non-empty and $\rho(G)$ is obtained from $\rho(G_1)$ and $\rho(G_2)$
by Klein combination.  If $\{\Lambda(\rho_j(G_m))\}$ converges to
$\Lambda(\rho(G_m))$ for both $m$, then $\{\rho_j\}$ converges
strongly to $\rho$ and $\{\Lambda(\rho_j(G))\}$ converges to
$\Lambda(\rho(G))$. 
\end{proposition}

\end{proof}

\section{Generalizations}

One may use Selberg's lemma (see Selberg \cite{selberg}) and results
of J\o rgensen and Marden to remove the assumption that $G$ is
torsion-free in the statement of Theorem \ref{strong-limit}.

\begin{theorem}{torsion}{}{} Let $G$ be a finitely generated group and
let $\{\rho_j\}\subset\calD(G)$ be a type-preserving sequence
converging algebraically to $\rho\in\calD(G)$.  If $\Omega(\rho(G))$
is non-empty, then $\{\rho_j\}$ converges strongly to $\rho$ and
$\{\Lambda(\rho_j(G))\}$ converges to $\Lambda(\rho(G))$.
\end{theorem}

\begin{proof}{} Selberg's lemma guarantees that there exists a finite
index, torsion-free subgroup $H$ of $G$.  Theorem \ref{strong-limit}
implies that $\{\rho_j|_H\} $ converges strongly to $\rho|_H$ 
and that $\{\Lambda(\rho_j(H))\}$ converges to $\Lambda(\rho(H))$.
Since $\Lambda(\rho_j(G))=\Lambda(\rho_j(H))$ for all $j$ and
$\Lambda(\rho(G))=\Lambda(\rho(H))$, we see that
$\{\Lambda(\rho_j(G))\}$ converges to $\Lambda(\rho(G))$.
If $G$ is not virtually abelian, then Proposition \ref{limcon}
implies that $\{\rho_j\}$ converges strongly to $\rho$.
If $G$ is virtually abelian, we may prove that $\{\rho_j\}$
converges strongly to $\rho$ by extending the arguments given
in the proof of Lemma \ref{elementary}.
\end{proof}

We can similarly generalize Corollary \ref{group-condition} to obtain:

\begin{corollary}{torsion2}{}{} Let $G$ be a finitely generated group
and let $\{\rho_j\}\subset\calD(G)$ be a type-preserving sequence
converging algebraically to $\rho$.  If $G$ contains a finite index,
torsion-free subgroup which is not a (non-trivial) free product of (orientable)
surface groups and infinite cyclic groups, then $\{\rho_j\}$ converges
strongly to $\rho$ and $\{\Lambda(\rho_j(G))\}$ converges to
$\Lambda(\rho(G))$.
\end{corollary}

\end{document}